\documentclass[conference]{IEEEtran}

\newtheorem{theorem}{\bf  Theorem}
\newtheorem{lemma}{\bf  Lemma}

\newtheorem{remark}{ \sc Remark}

\IEEEoverridecommandlockouts
\usepackage{cite}
\usepackage{amsmath,amssymb,amsfonts}
\usepackage{algorithmic}
\usepackage{graphicx}
\usepackage{textcomp}
\usepackage{xcolor}
\def\BibTeX{{\rm B\kern-.05em{\sc i\kern-.025em b}\kern-.08em
    T\kern-.1667em\lower.7ex\hbox{E}\kern-.125emX}}
\begin{document}

\title{An algorithm for J-spectral factorization of certain matrix functions
	\thanks{The authors were supported in part by Faculty Research funding from the Division of Science and Mathematics, NYUAD. 
		The first author was partially supported by the Shota Rustaveli National Science Foundation of
		Georgia (Project No. FR-18-2499).
	}
}

\author{\IEEEauthorblockN{Lasha~Ephremidze \textsuperscript{1,2} }
	\IEEEauthorblockA{1. \textit{Department of Mathematical Analysis} \\
		\textit{ Razmadze Mathematical Institute}\\
		Tbilisi, Georgia \\
		le23@nyu.edu}
	\and
		\IEEEauthorblockN{Ilya~Spitkovsky\textsuperscript{2} }	
		\IEEEauthorblockA{2. \textit{Division of Science and Mathematics} \\
		\textit{ New York University Abu Dhabi}\\
		Abu Dhabi, UAE \\
	ims2@nyu.edu}
}
\maketitle

\begin{abstract}
The problems of matrix spectral factorization and $J$-spectral factorization appear to be important for practical use in many MIMO control systems. We propose a numerical algorithm for $J$–spectral factorization which extends Janashia–Lagvilava matrix spectral factorization method to the indefinite case. The algorithm can be applied to matrices which have constant signatures for all leading principle submatrices. A numerical example is presented for illustrative purposes.
\end{abstract}

\begin{IEEEkeywords}
Spectral factorization,	$J$-spectral factorization,  algorithms.
\end{IEEEkeywords}

\section{Introduction}
Spectral factorization plays a prominent role in a wide
range of fields in system theory and control engineering.
In the scalar case, which arises in systems with single
input and single output, the factorization problem is relatively
easy and several classical methods exist to perform
this task (see a survey paper \cite{SayKai}). The matrix spectral
factorization, which arises in multi-dimensional systems,
is significantly more difficult. Following Wiener’s original
efforts \cite{Wie58}, dozens of papers addressed the development
of appropriate algorithms. None of the above methods can be implemented directly to solve the $J$-spectral factorization.

The Janashia–Lagvilava method is a relatively new algorithm for matrix spectral factorization \cite{JL99}, \cite{IEEE} which proved to be rather effective \cite{IEEE2018}. To describe this method of $r\times r$ matrix spectral factorization in a few words, one can say that it first performs a lower-upper triangular factorization with causal entries on the diagonal and then carries out an approximate spectral factorization of principle $m\times m$ submatrices step-by-step, $m=2,3,\ldots,r$. The decisive role in the latter process is played by unitary matrix functions
of certain structure, which eliminates many technical difficulties connected with computation. 

In the present paper, we extend Janashia-Lagvilava method to $J$-spectral factorization case, by using appropriately chosen $J$-unitary matrix functions instead of aforementioned unitary matrices. So far, the method can be used for matrices which have constant signatures for all leading principle submatrices, however, we hope to remove this restriction in the future work. Furthermore, the method has a potential of identifying a simple necessary and sufficient condition for the existence of $J$-spectral factorization and of being further extended towards the factorization of a wider class of Hermitian  matrices. 

Performed  numerical simulations confirm that  the proposed  algorithm,  whenever applicable, is as  effective  as  the  existing  matrix  spectral  factorization  algorithm.  On several occasions,  the  algorithm  can  also deal with the  so  called  singular cases,  where  the  zeros  of  the  determinant  occur  on  the  boundary.  Like the   Janashia–Lagvilava  method, the algorithm can  be  used  to  $J$-factorize  non–rational  matrices  as  well. 

\section{Formulation of the problem}

Let
\begin{equation}\label{1}
S(z)=\begin{pmatrix} s_{11}(z)& s_{12}(z)& \cdots&s_{1r}(z)\\
s_{21}(z)& s_{22}(z)& \cdots&s_{2r}(z)\\
\vdots&\vdots&\vdots&\vdots\\s_{r1}(z)& s_{r2}(z)&
\cdots&s_{rr}(z)\end{pmatrix},
\end{equation}
$z\in\mathbb{T}:=\{z\in\mathbb{C}:|z|=1\}$, be a Hermitian  $r\times r$ matrix function of constant signature, i.e. $S(z)=S^*(z)$ and the number of positive and negative eigenvalues of $S(z)$ are the constants $p$ and $q$, with $p+q=r$, for a.a. $z\in\mathbb{T}$.

$J$-spectral factorization of $S$ is by definition  the representation
\begin{equation}\label{2}
S(z)=S_+(z)\,J\, S_+^*(z),
\end{equation}
where $S_+$ can be extended to a stable analytic function inside $\mathbb{T}$, the matrix function $S_+^*$ is the Hermitian conjugate of $S$, and $J=(I_p\,,\;-I_q)$ is the diagonal matrix with $p$ ones and $q$ negative ones on the diagonal. We do not specify the classes to which $S$ and $S_+$ belong. For simplicity, one can assume that they are (Laurent) matrix polynomials. 

The necessity of factorization \eqref{2} arises in $\mathcal{H}_\infty$ control \cite{Fr87}, \cite{Ki97} and its solution is much more involved than the (standard) spectral factorization of positive definite matrix functions (when $p=r$ and $q=0$). Various algorithms for $J$-spectral factorization appear in the literature\cite{Se94}, \cite{St04} mostly for rational matrices.

Below, we present a new algorithm of $J$-spectral factorization which is an extension of Janashia-Lagvilava matrix spectral factorization method. Similarly to this method, we first perform a lower-upper triangular $J$-factorization of \eqref{1} with analytic entries on the diagonal.
This can be achieved only in the case where all the leading principal minors of $S$ have constant signs almost everywhere on $\mathbb{T}$, therefore, we impose this restriction on \eqref{1}.
 Then we recursively $J$-factorize leading principle $m\times m$ submatrices of $S$, $m=2,3,\ldots, r$.

\section{Notation}

For any set $\mathcal{S}$, we denote by $\mathcal{S}^{m\times n}$ the set of $m\times n$ matrices with entries from $\mathcal{S}$.

For a matrix $M\in\mathbb{C}^{r\times r}$ we use the standard notation $M^T$ and $M^*:=\overline{M}^T$ for the transpose and the Hermitian conjugate of $M$. The leading principle $m\times m$ submatrix of $M$, $m\leq r$, is denoted by $[M]_{m\times m}$. The same notation is used for matrix functions as well.

The letter $J$ always denotes a signature, i.e. a  square  diagonal  matrix  with  entries  $\pm1$  on the diagonal. The sizes  and  entries  of  $J$  may  vary  on different  occasions.  We  say  that  a  Hermitian   matrix 
$A=A^*\in\mathbb{C}^{m\times m}$  has  the  signature 
$J=(I_p\,,\;-I_q)$  if $A$ has $p$ positive and $q$ negative eigenvalues.

For a fixed signature matrix $J$, the set of $J$-unitary matrices, $\mathcal{U}_J$, is a group. Furthermore, $U\in\mathcal{U}_J \Longrightarrow U^T\in \mathcal{U}_J$, since $AJB=J \Longrightarrow  BJA=J$.

The set of polynomials is denoted by $\mathcal{P}^+$, and the set of Laurent polynomials,
\begin{equation}\label{14.04}
P(z)=\sum\nolimits_{k=-n}^m p_kz^k,
\end{equation}
is denoted by $\mathcal{P}$. The set of Laurent polynomials of degree at most $N$ (i.e. $0\leq n,m\leq N$ in \eqref{14.04}) is denoted by $\mathcal{P}_N$, and 
$$\mathcal{P}_N^+=\mathcal{P}_N \cap\mathcal{P}^+.$$

For Laurent polynomial \eqref{14.04}, let 
$$
\widetilde{P}(z)=\sum\nolimits_{k=-n}^m \overline{p_k}z^{-k}.
$$
Suppose also $\mathcal{P}_N^-:=\{P:\widetilde{P}\in\mathcal{P}_N^+\}$. Obviously, $\mathcal{P}_N^-\cap \mathcal{P}_N^+$ consists  of constant functions only.

A matrix polynomial $\mathbf{U}\in\mathcal{P}^{m\times m}$ is called $J$-unitary if $\mathbf{U}(z)$ is $J$-unitary for every $z\in\mathbb{T}$.

The $k$th Fourier coefficient of an integrable function $f\in
L_1(\mathbb{T})$ is denoted by $c_k\{f\}$.
If a function $f$ is square integrable, $f\in L_2=L_2(\mathbb{T})$, then
$$
f(z)=\sum\nolimits_{k=-\infty}^\infty c_k\{f\}z^k \;\;\text{ for a.a. }z\in\mathbb{T},
$$
and $\|f\|_2=2\pi \sum _{k=-\infty}^\infty |c_k\{f\}|^2$.

An integrable function $f$ is called analytic or causal if its Fourier expansion has the form
$$
f\sim \sum\nolimits_{k=0}^\infty c_k\{f\}z^k. 
$$
It is called stable if $f(z)\not=0$ for each $z$ with $|z|<1$, and it is called optimal if (see, e.g., \cite[Th. 17.17]{Rud}
$$
\log |f(0)|=\frac{1}{2\pi} \int_0^{2\pi} \log |f(e^{it})|\,dt.
$$

For a positive integrable function $f$ defined on $\mathbb{T}$, which satisfies the Paley-Wiener condition
$$
\log f\in L_1,
$$
there exists a unique (up to a constant multiple with absolute value 1) causal, stable, and optimal function $f^+$ such that
$$
f(z)= f^+(z)\overline{ f^+(z)}=|f^+(z)|^2  \;\text{ for a.a. }z\in\mathbb{T}.
$$
Such a function $f^+$ is called the (canonical) scalar spectral factor of $f$ and it can be given explicitly by the formula
$$
f^+(z)= \sqrt{f(z)}\exp\left(\frac12 i \mathcal{C}\big(\log f\big)(z)\right),
$$ 
where $\mathcal{C}$ stands for the harmonic conjugate of $f$:
$$
\mathcal{C} (f)(z)=\frac{1}{2\pi}(P) \int_0^{2\pi}f(e^{it}) \cot\frac{t-\tau}{2}\,d\tau,\;\;\;z=e^{it}.
$$	
This formula is the core of existing Exp-Log algorithm for scalar spectral factorization. It is the claim of well-known Fej\'er-Riesz lemma that if, in addition, $f\in\mathcal{P}_N$, then $f^+\in\mathcal{P}_N^+$. In Section V, we use the special  notation 
\begin{equation}\label{ssf}
f^+=\sqrt[+]{f}
\end{equation}
for the scalar spectral factor.

Finally, $\delta_{ij}$ stands for the  Kronecker delta, i.e. $\delta_{ij}=1$ if $i=j$ and $\delta_{ij}=0$ otherwise.

\section{The main observation}

In this section we generalize the main theorem of Janashia-Lagvilava method for $J$-unitary matrices. 

\begin{theorem} {\rm (cf. \cite[Th. 1]{IEEE})}  Let $F$ be an $m\times m$ matrix function of the form 	
	\begin{equation}\label{IE10}
	F=\begin{pmatrix}1&0&\cdots&0&0\\
	0&1&\cdots&0&0\\
	\vdots&\vdots&\vdots&\vdots&\vdots\\
	0&0&\cdots&1&0\\
	\zeta^-_{1}&\zeta^-_{2}&\cdots&\zeta^-_{m-1}&f^+
	\end{pmatrix},
	\end{equation}
	where
	\begin{equation}\label{IE11}
	\zeta^-_j\in \mathcal{P}_N^-,\;j=1,2,\ldots, m-1;\;
	f^{+}\in\mathcal{P}_N^+,\;f^+(0)\not=0,
	\end{equation}
	for some positive integer $N$, and let $J$ be an arbitrary signature. Then (almost surely) there exists a $J$-unitary matrix function $U$  of
	the form
	\begin{equation}\label{IE12}
	U=\begin{pmatrix}u_{11}&u_{12}&\cdots&u_{1m}\\
	u_{21}&u_{22}&\cdots&u_{2m}\\
	\vdots&\vdots&\vdots&\vdots\\
	u_{m-1,1}&u_{m-1,2}&\cdots&u_{m-1,m}\\[3mm]
	\widetilde{u_{m1}}&\widetilde{u_{m2}}&\cdots&\widetilde{u_{mm}}\\
	\end{pmatrix},
	\end{equation}
	where
	\begin{equation}
	u_{ij}\in \mathcal{P}_N^+,\;\;i,j=1,2,\ldots,m,
	\end{equation}
	with constant determinant, such that
	\begin{equation}\label{IE14}
	F U\in (\mathcal{P}^+_N)^{m\times m}.
	\end{equation}	
\end{theorem}

\begin{remark}
	A sketch  of  the  proof  below  indicates  the  isolated  cases  where  the  theorem  fails  to  hold. This is the sense in which we use the term ``almost surely”.    Whenever the solution exists, it is constructed explicitly.
\end{remark}

\smallskip

The  proof  follows  literally  the  proof  of  Theorem  1  in  \cite{IEEE}.  We  need  only  to  change  signs  of  some  expressions  accordingly.  By  this  way,  we  naturally  arrive  at  $J$-unitary  matrix  functions  instead  of  unitary  ones.
Indeed, for given functions $\zeta^-_j$, $j=1,2,\ldots,m-1$, 
$f^+$ satisfying \eqref{IE11}, and the signature $J={\rm diag}(J_1,J_2,\ldots,J_{m-1},1)$, we consider the following system of $m$
conditions (cf. (15) in \cite{IEEE})
\begin{equation}\label{IE15}
\begin{cases} \zeta^-_1x_m-J_1\cdot f^+\widetilde{x_1}\in \mathcal{P}^+,\\
\zeta^-_2x_m-J_2\cdot f^+\widetilde{x_2}\in \mathcal{P}^+,\\
\cdot\hskip+1cm \cdot\hskip+1cm \cdot\\
\zeta^-_{m-1}x_m-J_{m-1}\cdot f^+\widetilde{x_{m-1}}\in \mathcal{P}^+,\\
\zeta^-_1x_1+\zeta^-_2x_2+\ldots+\zeta^-_{m-1}x_{m-1}
+f^+\widetilde{x_m}\in \mathcal{P}^+,
\end{cases}
\end{equation}
where $\big(x_1,x_2,\ldots,x_m\big)^T\in (\mathcal{P}^+_N)^{m\times 1} $ is the unknown vector function. We say that a vector function
\begin{equation}
\mathbf{u}=\big(u_1,u_2,\ldots,u_m\big)^T\in
(\mathcal{P}^+_N)^{m\times 1}
\end{equation}
is a solution of \eqref{IE15} if and only if all the conditions in \eqref{IE15} are satisfied whenever $x_i=u_k$, $i=1,2,\ldots,m$.

We make essential use of the following
\begin{lemma}
	Let \eqref{IE11}  hold and  let
	\begin{gather*}
	\mathbf{u}=\big(u_1,u_2,\ldots,u_m\big)^T\in
	(\mathcal{P}_N^+)^{m\times 1} \\
	\mathbf{v}=\big(v_1,v_2,\ldots,v_m\big)^T\in
	(\mathcal{P}_N^+)^{m\times 1}
	\end{gather*}
	be two $($possibly identical$)$ solutions of the system \eqref{IE15}.
	Then
	\begin{equation}\label{IE19}
	\sum_{k=1}^{m-1}J_ku_k\widetilde{v_k}+\widetilde{u_m}v_m=\operatorname{const}.
	\end{equation}
\end{lemma}
{\em Proof:}
Substituting the functions $v$ in the first $m-1$ conditions and
the functions $u$ in the last condition of \eqref{IE15}, and then
multiplying the first $m-1$ conditions by $u$ and the last
condition by $v_m$,  we get
$$
\begin{cases} \zeta^-_1v_mu_1-J_1\cdot f^+\widetilde{v_1}u_1\in \mathcal{P}^+,\\
\zeta^-_2v_mu_2-J_2\cdot f^+\widetilde{v_2}u_2\in \mathcal{P}^+,\\
\cdot\hskip+1cm \cdot\hskip+1cm \cdot\\
\zeta^-_{m-1}v_mu_{m-1}-J_{m-1}\cdot f^+\widetilde{v_{m-1}}u_{m-1}\in \mathcal{P}^+,\\
\zeta^-_1u_1v_m+\zeta^-_2u_2v_m+\ldots+\zeta^-_{m-1}u_{m-1}v_m
+f^+\widetilde{u_m}v_m\in \mathcal{P}^+.
\end{cases}
$$
Subtracting the first $m-1$ conditions from the last condition in
the latter system, we get
\begin{equation}\label{IE20}
f^+\left(\sum_{k=1}^{m-1} J_k
u_k\widetilde{v_k}+\widetilde{u_m}v_m\right)\in
\mathcal{P}^+.
\end{equation}
Since the second multiple in \eqref{IE20} belongs to $\mathcal{P}_N$, 
taking into account the last condition in \eqref{IE11}, we get
\begin{equation*}
\sum_{k=1}^{m-1}J_ku_k\widetilde{v_k}+\widetilde{u_m}v_m\in
\mathcal{P}_N^+.
\end{equation*}
We can interchange the roles of $u$ and $v$ in the above
discussion to get in a similar manner that
\begin{equation*}
\sum_{k=1}^{m-1}J_k v_k\widetilde{u_k}+\widetilde{v_m}u_m\in
\mathcal{P}_N^+.
\end{equation*}
Consequently, the function in
\eqref{IE19} belongs to $\mathcal{P}_N^+\cap \mathcal{P}_N^-$, which
implies	\eqref{IE19}.  \hfill $\blacksquare$

The proof of Theorem 1 proceeds as follows. We search for a
nontrivial polynomial solution
\begin{equation}
\mathbf{x}=\big(x_1,x_2,\ldots,x_m\big)^T\in
(\mathcal{P}_N^+)^{m\times 1}
\end{equation}
of the system \eqref{IE15}, where
\begin{equation}\label{IE24}
x_i(z)=\sum_{n=0}^N a_{in} z^n,\;\;\;i=1,2,\ldots, m,
\end{equation}
and explicitly determine the coefficients $a_{in}$. We will find
such $m$ linearly independent  solutions of \eqref{IE15} which appear to be $m$ different columns of \eqref{IE12}

Equating all the  Fourier coefficients with non-positive indices
of the functions in the left-hand side of \eqref{IE15} to zero, except the
$0$th coefficient of the $j$th function which we set equal to $1$, we
get the following system of algebraic equations in the block
matrix form which we denote by $\mathbb{S}_j$:

\begin{equation}\label{IE25}
\mathbb{S}_j:=
\begin{cases}\Gamma_1 X_m-J_1D\overline{X_1}={\bf 0}, \\
\Gamma_2 X_m-J_2D\overline{X_2}={\bf 0}, \\
\;\;\;\;\;   \;\;\;\;\;    \\
\Gamma_j X_m-J_{j}D\overline{X_j}={\bf 1}, \\
\;\;\;\;\;   \;\;\;\;\;    \\
\Gamma_{m-1} X_m-J_{m-1}D\overline{X_{m-1}}={\bf 0}, \\
\Gamma_1 X_1+\Gamma_2 X_2+\ldots+\Gamma_{m-1}
X_{m-1}+D\overline{X_m}={\bf 0}\;. \end{cases}
\end{equation}
Here the following matrix notation is used:
\begin{gather*}
D=\begin{pmatrix}d_0&d_1&d_2&\cdots&d_{N-1}&d_N\\
0&d_0&d_1&\cdots&d_{N-2}&d_{N-1}\\
0&0&d_0&\cdots&d_{N-3}&d_{N-2}\\
\cdot&\cdot&\cdot&\cdots&\cdot&\cdot\\
0&0&0&\cdots&0&d_0\end{pmatrix},\;\;\\
\Gamma_i=
\begin{pmatrix}\gamma_{i0}&\gamma_{i1}&\gamma_{i2}
&\cdots&\gamma_{i,N-1}&\gamma_{iN}\\
\gamma_{i1}&\gamma_{i2}&\gamma_{i3}&\cdots&\gamma_{iN}&0\\
\gamma_{i2}&\gamma_{i3}&\gamma_{i4}&\cdots&0&0\\
\cdot&\cdot&\cdot&\cdots&\cdot&\cdot\\
\gamma_{iN}&0&0&\cdots&0&0\end{pmatrix},
\end{gather*}
$i=1,2,\ldots,m-1$, where
$$
f^+(z)=\sum_{n=0}^N d_n z^n \;\text{ and }\; \zeta^-_i(z)=
\sum_{n=0}^N\gamma_{in}z^{-n};
$$
\begin{equation*}
{\bf 0}=(0,0,\ldots,0)^T \text{ and } {\bf
	1}=(1,0,0,\ldots,0)^T\in \mathbb{C}^{N+1}.
\end{equation*}
The column vectors
\begin{equation*}
X_i=(a_{i0},a_{i1},\ldots,a_{iN})^T,\;\;i=1,2,\ldots,m,
\end{equation*}
(see \eqref{IE24}) are the unknowns.

Since $d_0=f^+(0)\not=0$ (see \eqref{IE11}), the matrix $D$ is invertible. Hence, determining $X_i$, $i=1,2,\ldots,m-1$, from the first $m-1$
equations of \eqref{IE25},
\begin{equation}\label{IE30}
X_i=J_i\left(\overline{D^{-1}}\;\overline{\Gamma_i}\;\overline{X_m}-\delta_{ij}\overline{D^{-1}}\;{\bf
	1}\right),
\end{equation}
$i=1,2,\ldots,m-1$, and then substituting them in the last equation of \eqref{IE25}, we get
\begin{gather*}
J_1\Gamma_1\,\overline{D^{-1}}\;\overline{\Gamma_1}\;\overline{X_m}+J_2\Gamma_2\,\overline{D^{-1}}\;\overline{\Gamma_2}\;
\overline{X_m}+\cdots\\
+J_{m-1}\Gamma_{m-1}\,\overline{D^{-1}}\;\overline{\Gamma_{m-1}}\;\overline{X_m}+D\;\overline{X_m}=J_j\Gamma_j\,\overline{D^{-1}}\,{\bf 1}
\end{gather*}
(it is  assumed that the right-hand
side is equal to ${\bf 1}$ when $j=m$) or, equivalently,
\begin{gather}
(J_1\Theta_1\,{\Theta_1^*}+J_2\Theta_2\,{\Theta_2^*}
+\!\ldots\!+J_{m-1}\Theta_{m-1}\,{\Theta_{m-1}^*}+I_{N+1})\,\overline{X_m}\notag
\\
=J_jD^{-1}\,\Gamma_j\,\overline{D^{-1}}\,{\bf 1}, \label{IE31}
\end{gather}
where
\begin{equation*}
\Theta_i=D^{-1}\,\Gamma_i\,,\;\;i=1,2,\ldots,m-1
\end{equation*}
(we wrote $\Theta^*$ instead of $\overline{\Theta}$ because $\Theta^T=\Theta$). 

For each $j=1,2,\ldots,m$, \eqref{IE31} is a linear algebraic system of
$N+1$ equations with $(N+1)$ unknowns. This system \eqref{IE31} and consequently \eqref{IE25} has the unique solution for each $j=1,2,\ldots,m$ if and only if 
\begin{equation}\label{dtn0}
\det(\Delta)\not=0,\;\text{ where }\; \Delta=\sum\nolimits_{k=1}^{m-1}J_k\Theta\Theta^*+I_{N+1}\,.
\end{equation}

\begin{remark}
Unlike  the  spectral   factorization,  where   $\Delta$  is  always positive  definite  and  \eqref{dtn0}  holds,  there are isolated  indefinite  cases  where  \eqref{dtn0}  does  not  hold.  However,  we  can  assume  that  \eqref{dtn0}   holds  (see Remark 1)  and   proceed  with  solution  of  \eqref{IE15}.
\end{remark}

\begin{remark}
	As in the spectral factorization case (see \cite[Appendix]{IEEE}) the matrix $\Delta$ has a displacement structure of rank $m$ with respect to $Z$, where $Z$  is the upper triangular $(N+1)\times(N+1)$ matrix with 1’s on the first superdiagonal and 0’s elsewhere (i.e., a Jordan block with eigenvalue 0). Namely,
$$
R_Z\Delta:=\Delta-Z\Delta Z^*=AJA^*,
$$
where $A$ is the $(N+1)\times m$ matrix which has $i$-th column equal to the first column of $\Theta_i$, $i=1,2,\ldots,m-1$, and the last column is equal to $(0,0,\ldots,0,1)\in\mathbb{C}^{N+1}$. Consequently, the triangular factorization of $\Delta$ can be performed in $O(mN^2)$ operations instead of the traditional $O(N^3)$ ones, as it is described in \cite[Appendix F.1]{Kai99}. This substantially reduces the amount of operations if $N\gg m$.
\end{remark}

Finding the matrix vector $\overline{X_m}$ from \eqref{IE31} and then
determining $X_1,X_2,\ldots,X_{m-1}$ from \eqref{IE30}, we get the unique
solution of $\mathbb{S}_j$. To indicate its dependence on $j$, we
denote the solution of $\mathbb{S}_j$ by
$(X_1^j,X_2^j,\ldots,X_{m-1}^j,X_m^j)$,
\begin{equation}\label{IE35}
X_i^j:=(a_{i0}^j,a_{i1}^j,\ldots, a_{iN}^j)^T,
\;\;\;i=1,2,\ldots,m,
\end{equation}
so that if we construct a matrix function $V$,
\begin{equation}\label{IE36}
V=\begin{pmatrix}v_{11}&v_{12}&\cdots&v_{1m}\\
v_{21}&v_{22}&\cdots&v_{2m}\\
\vdots&\vdots&\vdots&\vdots\\
v_{m-1,1}&v_{m-1,2}&\cdots&v_{m-1,m}\\[3mm]
\widetilde{v_{m1}}&\widetilde{v_{m2}}&\cdots&\widetilde{v_{mm}}\\
\end{pmatrix},
\end{equation}
by letting (see \eqref{IE35})
\begin{equation}
v_{ij}(z)=\sum_{n=0}^N a_{in}^j z^n, \;\;\;1\leq i,j\leq m,
\end{equation}
then columns of \eqref{IE36} are  solutions of the system \eqref{IE35}.  Hence, because
of the last equation in (15),
\begin{equation*}
FV\in (\mathcal{P}^{+}_N)^{m\times m}
\end{equation*}
and, by virtue of Lemma 1,
\begin{equation}\label{Const2}
V(z)\,J\,V^*(z)=C,
\end{equation}
where $C$ is a constant Hermitian  matrix with signature $J$. It can be also proved that (see \cite[p. 2322, II]{IEEE}) that
$$
\det V(z)={\rm const}.
$$
Decomposing the matrix $C$ as 
\begin{equation}\label{CJC}
C=C_0J\,C_0^*,\;\text{ where }C_0=V(1),
\end{equation}
equations \eqref{Const2} and \eqref{CJC} imply
$$
C_0^{-1}V(z)\,J\,(C_0^{-1}V(z))^*=J.
$$
Hence, 
$$U=C_0^{-1}V$$
is the required $J$-unitary matrix and it can be numerically computed by using the above equations. \hfill  $\blacksquare$

\section{Description of the algorithm}

In this section we provide computational procedures for $J$-factorization of \eqref{1} which are similar to corresponding procedures presented in \cite{IEEE}.

{\bf Procedure 1.} First we perform the lower-upper triangular
$J$-factorization of $S$:

\begin{equation} \label{IE54}
S(z)=M(z)\,J\,M^*(z).
\end{equation}
Here 
$$
M=\begin{pmatrix}f^+_1&0&\cdots&0&0\\
\xi_{21}&f^+_2&\cdots&0&0\\
\vdots&\vdots&\vdots&\vdots&\vdots\\
\xi_{r-1,1}&\xi_{r-1,2}&\cdots&f^+_{r-1}&0\\
\xi_{r1}&\xi_{r2}&\cdots&\xi_{r,r-1}&f^+_r
\end{pmatrix},
$$
where  $f^+_m$,  $m=1,2,\ldots,r$,  are stable  analytic functions (we  also  assume  that  all  entries  are  square  integrable).  Such  factorization  can  always  be  achieved  under  the   restriction  that  
$\det [S]_{m\times m}$ has constant sign almost everywhere on $\mathbb{T}$ for each $m=1,2,\ldots,r$.
This  happens,  for example,  if  all  principle  minors  are  non-singular everywhere on $\mathbb{T}$, however,  this  condition  is not  necessary.  We  can  apply  the  similar  recursive  formulas  as  for  usual  Cholesky  factorization:   $f_1^+=
\sqrt[+]{J_1s_{11}} $,  $\xi_{i1}=J_1s_{i1}/\overline{f_1^+}$, $i=2,3,\ldots,r$;
\begin{gather*}
{f_j^+}=\sqrt[+]{J_j\left(s_{jj}-\sum\nolimits_{k=1}^{j-1}J_k\xi_{jk}
	\overline{\xi_{jk}}\right)},\; j=2,3,\ldots,r;
\\
\xi_{ij}=J_j\left(s_{ij}-\sum\nolimits_{k=1}^{j-1}J_k\xi_{ik}
\overline{\xi_{jk}}\right)/\overline{f_j^+},
\end{gather*}
$j=2,3,\ldots,r-1$, $i=j+1,j+2,\ldots,r$, assuming  that  $\sqrt[+]{\cdot}$  performs  the  scalar  spectral  factorization (see \eqref{ssf}).
In  actual computations, one  can  perform  factorization  \eqref {IE54}  pointwise  in  frequency  domain  for  selected  values  of  $z\in\mathbb{T}$.

\smallskip

{\bf Procedure 2.}  We approximate $M$ in $L_2$ keeping only a
finite number of coefficients with negative indices in the Fourier
expansions of the entries of $M$. For the convenience of
computations, we take a different number of these coefficients for
different entries below the main diagonal. Namely, for a large
positive integer $N$, let
\begin{equation}
M_N=\begin{pmatrix}f^+_1&0&\cdots&0&0\\[1mm]
\xi^{[N]}_{21}&f^+_2&\cdots&0&0\\[1mm]
\vdots&\vdots&\vdots&\vdots&\vdots\\[1mm]
\xi^{[(r-2)N]}_{r-1,1}&\xi^{[(r-3)N]}_{r-1,2}&\cdots&f^+_{r-1}&0\\[1mm]
\xi^{[(r-1)N]}_{r1}&\xi^{[(r-2)N]}_{r2}&\cdots&\xi^{[N]}_{r,r-1}&f^+_r
\end{pmatrix}
\end{equation}
where $\xi^{[{N}]}_{ij}(z)=\sum_{n=-{N}}^\infty
c_n\{\xi_{ij}\}z^n$, $2\leq i\leq r$, $1\leq j<r$. Let
$$
S_N(z)=M_N(z)\,J\,M_N^*(z).
$$

\smallskip

{\bf Procedure 3.} We compute explicitly $S_N^+$, a $J$-spectral
factor of $S_N$.
This is done recursively with respect to $m$. Namely, we represent $S_N^+$ as
$$
S_N^+=M_N\mathbf{U}_1\mathbf{U}_2\mathbf{U}_3\ldots
\mathbf{U}_r,
$$
where each $U_m$ is $J$-unitary and has the block matrix form
\begin{equation}\label{IE63}
\mathbf{U}_m(t)=\begin{pmatrix}U_{m}(t)&0\\0&I_{r-m}\end{pmatrix},
\end{equation}
$m=2,3,\ldots r$. Furthermore, each $[Q_m]_{m\times m}$ is $J$-spectral factor of $[S_N]_{m\times m}$
\begin{equation}\label{SJQ}
[S_N]_{m\times m}=[Q_m]_{m\times m}\,[J]_{m\times m} \,[Q_m]_{m\times m}^*,
\end{equation}
where 
$$
Q_m=M_N\mathbf{U}_1\mathbf{U}_2\mathbf{U}_3\ldots
\mathbf{U}_m.
$$
We take $ \mathbf{U}_1=I_r$ and then \eqref{SJQ} is valid for $m=1$. Assume that $\mathbf{U}_2(t)$,$\mathbf{U}_2(t)${}$\ldots$,$\mathbf{U}_{m-1}(t)$
have already been constructed so that \eqref{SJQ} holds when $m$ is replaced by $m-1$ and suppose the last row of 
$[Q_{m-1}]_{m\times m}$ is
$[\zeta^{m-1}_{1},\zeta^{m-1}_{2},\ldots,\zeta^{m-1}_{m-1},f_m^+]$. 
Then we construct the next $J$-unitary matrix \eqref{IE63} by performing the following operations: 

\begin{figure}[htbp]
	\centerline{\includegraphics[width=95mm,scale=1.2]{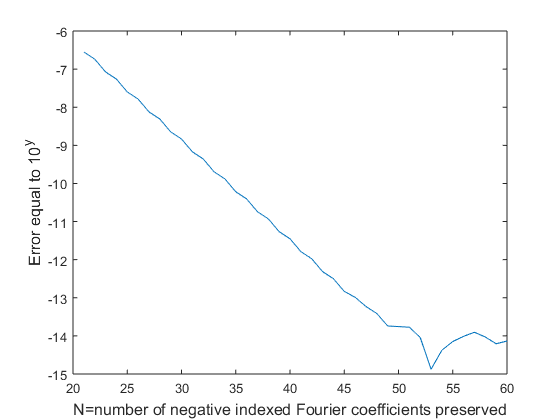}}
	\caption{Error in $J$-spectral factorization of matrix (29)}
	\label{fig}
\end{figure}

{\sc  Step 1.} Construct a matrix function $F(t)$ of the form \eqref{IE10}, where
$$\zeta^-_j(z)=\sum_{n=-(m-1)N}^0
c_n\big\{\zeta^{m-1}_j\big\}\,z^n, \;\;\;j=1,2,\ldots,m-1,
$$
and
$$f^+(z)=\sum\nolimits^{(m-1)N}_{n=0}c_n\{f_m^+\}\,z^n\,.$$

{\sc Step 2.} Using Theorem 1, construct $U$ of the form \eqref{IE12}, where
$u_{ij}\in\mathcal{P}^+_{(m-1)N}$, $1\leq i,j\leq m$, so that
\eqref{IE14} would hold.

{\sc Step 3.} Define $\mathbf{U}_m$ by the equation \eqref{IE63} where
$U_m=U$ is found in Step 2.

\section{Numerical Example}

To illustrate our approach, we present an approximate $J$-factorization of the following polynomial matrix function $S=$
\begin{equation}\label{14.01}
\begin{pmatrix}    
-8z^{-1}-19-8z & -39z^{-1}-73-28z\\
-28z^{-1}-73-39z & -137z^{-1}-286-137z
\end{pmatrix}.
\end{equation}
This matrix satisfies the conditions imposed on $S$ in order for the algorithm to be applicable, namely $s_{11}(z)$ and
$$
\det S(z)=4(z^{-2}-2+z^2)=4(z^{-2}-1)(1-z^2)
$$
are both negative for $z\in\mathbb{T}$. 
However, the matrix $S(z)$ is singular for $z=-1$ and $1$, which usually complicates the factorization process. The $J$-factorization of  \eqref{14.01} is known in advance due to the corresponding example of the singular matrix in \cite{IEEE2018}: 
$$
S(z)=S_+(z)\begin{pmatrix}    
-1 & 0\\ 0& 1
\end{pmatrix} S_+^*(z),
$$
where 
\begin{equation}\label{14.02}
S_+(z)=
\begin{pmatrix}    
4     +2z & 1\\  14 +   10z &  3  +   z
\end{pmatrix}.
\end{equation}
However, we follow the steps of the proposed algorithm to produce an approximate result.

The triangular $J$-factorization of $S$ has the form
$$
S(z)=M(z)
\begin{pmatrix}    
-1 & 0\\ 0& 1
\end{pmatrix}
M^*(z)
$$
where $M(z)=$
$$
\begin{pmatrix}    
3.824\ldots+z\cdot2.092\ldots  & 0\\[2mm]\dfrac{28z^{-1}+73+39z}{z^{-1}\cdot2.092\ldots+ 3.824\ldots} & \dfrac{1-z^2}{3.824\ldots+z\cdot2.092\ldots}
\end{pmatrix}
$$
with
$$
f_1^+(z):= 3.824\ldots+z\cdot2.092\ldots=\sqrt[+]{8z^{-1}+19+8z}
$$
and
$$
f_2^+(z):=(1-z^2)=\sqrt[+]{-z^{-2}+2-z^2}.
$$
We expand $\xi_{21}=-s_{21}/\widetilde{f^+_1}$ into Fourier series by the  division of polynomials and, for a positive integer $N$, approximate it by ``cutting the tail”:
$$
\xi_{21}(z)\approx\xi_{21}^{[N]}(z)=\sum\nolimits_{k=-N}^{\infty}c_k\{\xi_{21}\}z^k.
$$
Thus we get the approximation of $S$ by
$$
S_N=\begin{pmatrix}    
f_1^+ & 0\\ \xi_{21}^{[N]}& f_2^+
\end{pmatrix}
\begin{pmatrix}    
-1 & 0\\ 0& 1
\end{pmatrix}
\begin{pmatrix}    
f_1^+ & 0\\ \xi_{21}^{[N]}& f_2^+
\end{pmatrix}^*
$$
and we obtain its $J$-spectral factor $S_N^+$ by finding explicitly a $J$-unitary matrix $U=U_N$ as it is described in Section V:
$$
S_N^+=\begin{pmatrix}    
f_1^+ & 0\\ \xi_{21}^{[N]}& f_2^+
\end{pmatrix}\cdot U.
$$
The computation results coincide with the exact answer \eqref{14.02} within 16 digits (the Matlab double precision) for $N=53$.

\vfill\break

 A total computational time to achieve this accuracy is less than 0.02 sec (on a laptop with the  characteristics: Intel(R) Core(TM) i7 8650U CPU, 1.90 GHz, RAM 16.00 Gb). Fig. 1 shows how this accuracy increases  with increasing $N$.

\section*{Acknowledgment} 
Authors thank Professor Michael \v Sebek for bringing to their attention the importance of $J$-spectral factorization in Control Theory.


\def\cprime{$'$}

\end{document}